\documentclass[12pt]{article}
\usepackage{amsthm,amsfonts}
\usepackage{fullpage}

\newtheorem{theorem}{Theorem}

\title{A note on avoidable words in squarefree ternary words}
\author{Narad Rampersad \\
School of Computer Science \\
University of Waterloo \\
Waterloo, ON, N2L 3G1 \\
CANADA \\
{\tt nrampersad@math.uwaterloo.ca}}

\begin{document}
\date{\today}
\maketitle

\begin{abstract}
We completely characterize the words that can be avoided in
infinite squarefree ternary words.
\end{abstract}

\section{Introduction}
Let $\Sigma$ be a finite, non-empty set called an \emph{alphabet}.
We denote the set of all words of finite length over the alphabet
$\Sigma$ by $\Sigma^*$.  Let $\Sigma_k$ denote the alphabet
$\{0,1,\ldots,k-1\}$; \emph{e.g.}, $\Sigma_3=\{0,1,2\}$.

A map $h:\Sigma^* \rightarrow \Delta^*$ is called a \emph{morphism} if
$h$ satisfies $h(xy)=h(x)h(y)$ for all $x,y\in\Sigma^*$.  A morphism
may be defined simply by specifying its action on $\Sigma$.
A morphism $h:\Sigma^* \rightarrow \Sigma^*$ such that $h(a)=ax$
for some $a\in\Sigma$ is said to be \emph{prolongable on $a$}; we
may then repeatedly iterate $h$ to obtain the \emph{fixed point}
$h^\omega(a)=axh(x)h^2(x)h^3(x)\cdots$.

An \emph{square} is a word of the form $xx$, where $x\in\Sigma^*$.
A word $w'$ is called a subword of $w$ if $w$
can be written in the form $uw'v$ for some $u,v\in\Sigma^*$.
We say a word $w$ is \emph{squarefree} (or \emph{avoids squares})
if no subword of $w$ is an square.

It is easy to check that no binary word of length $\geq 4$ avoids
squares.  However, Thue \cite {Thu06} gave an example of a infinite
squarefree ternary word.  There are certain words that are avoidable
in infinite squarefree ternary words and others that are unavoidable;
\emph{e.g.}, the word 101 is avoidable, whereas the word 012 is not.
In the next section we characterize all words that can be avoided in
infinite squarefree ternary words.

\section{Results}
\begin{theorem}
Let $w$ be any infinite squarefree word over $\Sigma_3$.  Then $w$
contains at least one occurrence of each of the following words:
012, 021, 102, 120, 201, 210.
\end{theorem}

\begin{proof}
This can be verified by an exhaustive computer search.  It suffices
to check all 34422 squarefree words of length 30 over $\Sigma_3$.
\end{proof}

\begin{theorem}
\label{abca}
Let $a$, $b$, and $c$ be distinct letters of $\Sigma_3$.  Then there
exists an infinite squarefree word over $\Sigma_3$ that
contains no occurrences of each of the words $abca$ and $acba$.
\end{theorem}

\begin{proof}
It is easy to see that $(abca,acba) \in \lbrace(0120,0210),
(1021,1201),(2012,2102)\rbrace$.  Hence it suffices to show that
there exists an infinite squarefree word over $\Sigma_3$ that avoids
0120 and 0210, as we may simply rename $a$, $b$, and $c$ to get the
desired avoidance.  Consider the morphism $h$ defined as follows:
\begin{eqnarray*}
0 & \rightarrow & 12  \\
1 & \rightarrow & 102 \\
2 & \rightarrow & 0
\end{eqnarray*}
Then the fixed point $h^\omega(0)$ is squarefree and avoids 101 and 202.
The only way to obtain 0120 from the morphism $h$ is $h(202) = 0120$, but
$h^\omega(0)$ avoids 202.  Similarly, the only way to obtain 0210 is
$1h(11)2 = 102102$, but $h^\omega(0)$ avoids the square 11.  The
result now follows.
\end{proof}

\begin{theorem}
Let $x$ be any word over $\Sigma_3$ such that $|x| \geq 4$.  Then there
exists an infinite squarefree word over $\Sigma_3$ that contains no
occurrences of $x$.
\end{theorem}

\begin{proof}
It suffices to prove the theorem for $|x| = 4$.
Consider the set $\mathcal{A}$ of all squarefree words of length
4 over $\Sigma_3$.  We have $\mathcal{A} = \mathcal{A}' \cup \mathcal{A}''$,
where $$\mathcal{A}' = \lbrace 0102,0121,0201,0212,1012,1020,1202,1210,
2010,2021,2101,2120 \rbrace$$ and $$\mathcal{A}'' = \lbrace 0120,0210,
1021,1201,2012,2102 \rbrace.$$  Note that all words in $\mathcal{A}'$
contain a subword of the form $aba$, where $a$ and $b$ are distinct
letters of $\Sigma_3$.  It is well known that for any such subword $aba$,
there exists an infinite squarefree word over $\Sigma_3$ that avoids $aba$.
Hence, for any word $x \in \mathcal{A}'$, there exists an infinite squarefree
word over $\Sigma_3$ that avoids $x$.

Now consider the set $\mathcal{A}''$.  Note that all words in $\mathcal{A}''$
are of the form $abca$, where $a$, $b$, and $c$ are distinct letters of
$\Sigma_3$.  Theorem~\ref{abca} implies that for any such word $abca$,
there exists an infinite squarefree word over $\Sigma_3$ that avoids
$abca$.  Hence, for any word $x \in \mathcal{A}$, there exists an infinite
squarefree word over $\Sigma_3$ that avoids $x$.  The result now follows.
\end{proof}

\end{document}